\documentclass[global,final]{svjour}
\usepackage{amsmath,amsfonts,amssymb}
\usepackage[all]{xy}
\usepackage[numbers,square,sort&compress]{natbib} 
\usepackage{xspace}
\date{12 November, 2007}
\def\Z{\mathbb{Z}}
\def\simkappa{\!\!\sim_{\kappa}}
\def\Acyc{\mathsf{Acyc}}
\def\acyc{\alpha}
\def\card#1{|#1|}

\def\vset{\mathrm{v}}
\def\eset{\mathrm{e}}
\def\CycleSet{\operatorname{Cycle}}
\def\eg{e.g.\xspace}
\def\ie{i.e.\xspace}
\def\makeheadbox{}
\begin{document}
\title{On Enumeration of Conjugacy\\Classes of Coxeter Elements}
\author{
Matthew~Macauley\inst{1,3} \and
Henning~S.~Mortveit\inst{2,3}}
\institute{
Department of Mathematics, UCSB \and 
Department of Mathematics, Virginia Tech \and 
NDSSL, VBI, Virginia Tech, \email{\{macauley, henning.mortveit\}@vt.edu}
}
\maketitle

\medskip

\abstract{
  In this paper we study the equivalence relation on the set of
  acyclic orientations of a graph $Y$ that arises through
  source-to-sink conversions. This source-to-sink conversion encodes,
  e.g. conjugation of Coxeter elements of a Coxeter group. We give a
  direct proof of a recursion for the number of equivalence classes of
  this relation for an arbitrary graph $Y$ using edge deletion and
  edge contraction of non-bridge edges. We conclude by showing how
  this result may also be obtained through an evaluation of the Tutte
  polynomial as $T(Y,1,0)$, and we provide bijections to two other
  classes of acyclic orientations that are known to be counted in the
  same way. A transversal of the set of equivalence classes is given.
}


\section{Introduction}
\label{sec:intro}

The equivalence relation on the set of acyclic orientations of a graph
$Y$ that arises from iteratively changing sources into sinks appears
in many areas of mathematics. For example, in the context of Coxeter
groups the source-to-sink operation encodes conjugation of Coxeter
elements~\cite{Shi:01}, although in general, these conjugacy classes
are not fully understood. Additionally, it is closely related to the
reflection functor in the representation theory of
quivers~\cite{Marsh:03}. It has also been studied in the context of
the chip-firing game of Bj\"{o}rner, Lov\'{a}sz, and
Shor~\cite{Bjorner:91}. Moreover, it arises in the characterization of
cycle equivalence for a class of discrete dynamical
systems~\cite{Macauley:07b}, which was the original motivation for
this work.

In~\cite{Shi:01}, the number of equivalence classes $\kappa(Y)$ for a
graph $Y$ was determined for graphs that contain precisely one cycle.
The main result of this paper is a novel proof for $\kappa(Y)$ for arbitrary
graphs in the form of a recurrence relation involving the edge
deletion $Y_e'$ and edge contraction $Y''_e$ of a cycle-edge $e$ in
$Y$. It can be stated as follows.

\begin{theorem}
\label{thm:collapse}
  Let $e$ be a cycle-edge of $Y$. Then
  \begin{equation}
    \label{eq:recursion}
    \kappa(Y)=\kappa(Y_e')+\kappa(Y_e'')\;.
  \end{equation}
\end{theorem}

Our proof involves a careful consideration of what happens to the
$\kappa$-equivalence classes of $\Acyc(Y)$ as a cycle-edge $e$ is
deleted. This leads to the construction of the collapse graph of $Y$
and $e$, which has vertex set the $\kappa$-classes of $\Acyc(Y)$. We
show that there is a bijection from the set of connected components of
this graph to the set of $\kappa$-equivalence classes of
$\Acyc(Y_e')$. Moreover, we establish that there is a bijection from
the edge set of the collapse graph to the set of $\kappa$-equivalence
classes of $\Acyc(Y_e'')$. From this and the fact that the collapse
graph is a forest the recursion~\eqref{eq:recursion} follows.
Alternatively, the recursion can be derived through an observation
made by Vic Reiner, see~\cite[Remark 5.5]{Novik:02}: the number of
equivalence classes of linear orderings under the operations of $(i)$
transposition of successive, non-connected generators and $(ii)$
cyclic shifts is counted by~\eqref{eq:recursion}. The bijection
between Coxeter elements and acyclic orientations in~\cite{Shi:01}
provides the connection to our setting. Even though the connection of
this fact to the enumeration of conjugacy classes of Coxeter elements
is straightforward, this does not appear in the literature.
Our contribution is an independent and direct proof of this result by
examining the acyclic orientations of the Coxeter graph. Additionally,
our proof provides insight into the structure of the equivalence
classes. We believe that the techniques involved may be useful in
extending current results in Coxeter theory, in particular, some
from~\cite{Shi:01}. 
\medskip

Let $Y$ be a finite undirected graph with vertex set $\vset[Y] =
\{1,2,\dots,n\}$ and edge set $\eset[Y]$. An orientation of $Y$ is
represented by a map $O_Y \colon \eset[Y] \longrightarrow \vset[Y]
\times \vset[Y]$, and the graph $G(O_Y)$ is obtained from
$Y$ by orienting each edge as given by $O_Y$. We will use $O_Y$ and
$G(O_Y)$ interchangeably when no ambiguity can arise. An orientation
$O_Y$ is acyclic if $G(O_Y)$ has no directed cycles. The set of
acyclic orientations of $Y$ is denoted $\Acyc(Y)$, and we set
$\acyc(Y)=\card{\Acyc(Y)}$, which can be computed through the
well-known recursion relation
\begin{equation}
  \label{eq:acyclic}
  \acyc(Y) = \acyc(Y'_e)+\acyc(Y''_e)\;.
\end{equation}
As above, $Y'_e$ and $Y''_e$ are the graphs obtained from $Y$ by
deletion and contraction of a fixed edge $e$, respectively.  It is
known that there is a bijection between $\Acyc(Y)$ and the set of
Coxeter elements of Coxeter group whose Coxeter graph is
$Y$~\cite{Humphreys:92,Shi:97a}. There is also a bijection between
$\Acyc(Y)$ and the set of chambers of the graphic hyperplane
arrangement $\mathcal{H}(Y)$~\cite{Orlik:92}.

If $v$ is a source of an acyclic orientation $O_Y$ with degree $\ge
1$, then reversing the orientation of all the edges incident to $v$
maps $O_Y$ to a new orientation of $Y$, which is also acyclic. This is
called a \emph{source-to-sink} operation, or a \emph{click}. We define
the equivalence relation $\sim_\kappa$ on the set of acyclic
orientations for a fixed graph $Y$ by $O_Y\sim_\kappa O'_Y$ if there
is a sequence of source-to-sink operations that maps $O_Y$ to
$O'_Y$. Two such orientations are said to be \emph{click-equivalent},
or \emph{$\kappa$-equivalent}. We set $\kappa(Y) =
\card{\Acyc(Y)/\simkappa}$.

The set of linear orders on $\vset[Y]$ can be represented by the set
of permutations of $\vset[Y]$, which we denote as $S_Y$. We write
$[\pi]_Y$ for the set of linear orders compatible with the acyclic
orientation $O_Y$ induced by $\pi$. There is a bijection between
$\bigl\{ [\pi]_Y \mid \pi\in S_Y \bigr\}$ and $\Acyc(Y)$, see,
\eg~\cite{Reidys:98a}.
Let $\pi$ be the permutation representation of a linear order
compatible with $O_Y$. Note that mapping $\pi =
(\pi_1,\pi_2,\ldots,\pi_n)$ to $\pi' = (\pi_2,\ldots,\pi_n,\pi_1)$
corresponds to converting $\pi_1$ from a source to a sink in $O_Y$. In
general, two distinct acyclic orientations $O_Y$ and $O_Y'$ are
$\kappa$-equivalent if and only if there exists $\pi$ compatible with
$O_Y$ and $\pi'$ compatible with $O_Y'$ such that $\pi'$ can be
obtained from $\pi$ by $(i)$ cyclic shifts and $(ii)$ transpositions
of consecutive elements that are not connected in $Y$. 
For a given Coxeter group $W$ with generators $S=\{s_i\}_{i=1}^n$ and
Coxeter graph $Y$ there is a similar mapping from $S_Y$ into the set
of Coxeter elements $C(W)$, and a bijection from $C(W)$ to $\Acyc(Y)$,
see~\cite{Shi:97a}.  Thus, an acyclic orientation represents a unique
Coxeter element, and a source-to-sink operation corresponds to
conjugating that element by a particular generator. Therefore,
$\kappa(Y)$ is an upper bound for the number of conjugacy classes of
Coxeter elements in a Coxeter group whose Coxeter graph is $Y$, and
this bound is known to be sharp in certain cases~\cite{Shi:01}. A
simple induction argument shows that if $Y$ is a tree, then
$\kappa(Y)=1$, and thus all Coxeter elements in a finite Coxeter group
are conjugate~\cite{Humphreys:92}. In~\cite{Shi:01}, the author shows
that if $Y$ contains a single cycle of length $n$, then
$\kappa(Y)=n-1$.  This becomes a straightforward corollary of
Theorem~\ref{thm:collapse}. The recurrence of
Theorem~\ref{thm:collapse} appears in several areas of mathematics,
and corresponds to the evaluation of the Tutte polynomial at $(1,0)$,
which we describe in Section~\ref{sec:tutte}.


\section{Preliminary Results}
\label{sec:prelim}

We begin our study of $\kappa(Y)$ by making the following simple
observation recorded without proof.

\begin{proposition} 
\label{prop:dunion}
Let $Y$ be the disjoint union of undirected graphs $Y_1$ and
$Y_2$. Then 
\begin{equation}
\label{eq:k-dunion}
\kappa(Y)=\kappa(Y_1)\kappa(Y_2)\;.
\end{equation}
\end{proposition}

In light of this, we may assume that $Y$ is connected when computing
$\kappa(Y)$. An edge of $Y$ that is not contained in any simple cycle
of $Y$ is a \emph{bridge}, otherwise it is a \emph{cycle-edge}. The
graph obtained from $Y$ by deletion of all bridges is the \emph{cycle
graph} of $Y$ and it is denoted $\CycleSet(Y)$. Alternatively, an edge
$e$ of a connected graph $Y$ is a bridge if the deletion of $e$
disconnects $Y$. Bridges do not contribute to $\kappa(Y)$ as shown in
the following proposition.

\begin{proposition}
  \label{prop:non-cycle-edge}
  Let $Y$ be an undirected and graph, and let $e=\{v,w\}$ be a bridge
  of $Y$, connecting the disjoint subgraphs $Y_1$ and $Y_2$. Then one
  has the relation
\begin{equation}
\label{eq:k-bridge}
\kappa(Y) = \kappa(Y_1) \kappa(Y_2)\;.
\end{equation}
\end{proposition}

\begin{proof}
Each pair of acyclic orientations $O_{Y_1} \in \Acyc(Y_1)$ and
$O_{Y_2}\in\Acyc(Y_2)$ extends to exactly two acyclic orientations of
$Y$ by $O_Y=(O_{Y_1},(v,w),O_{Y_2})$ and $O_Y' =
(O_{Y_1},(w,v),O_{Y_2})$ defined in the obvious way. Clearly,
every acyclic orientation of $Y$ is also of one of these forms.
Moreover, any click sequence for $O_Y'$ that contains each vertex of
$Y_2$ exactly once and contains no vertices of $Y_1$ maps $O_Y'$ to
$O_Y$. Hence $O_Y$ and $O_Y'$ are click-equivalent. It follows that
$O_Y, O_Y'\in\Acyc(Y)$ are click-equivalent if and only their
corresponding acyclic orientations over $Y_1$ and $Y_2$ are
click-equivalent, and the equality~\eqref{eq:k-bridge} now follows
from Proposition~\ref{prop:dunion}.
\qed
\end{proof}

Proposition~\ref{prop:non-cycle-edge} gives us the immediate corollary.

\begin{corollary}
  \label{cor:pruning}
  For any undirected graph $Y$ we have $\kappa(Y) =
  \kappa(\CycleSet(Y))$. In particular, if $Y$ is a forest then
  $\kappa(Y) = 1$. 
\end{corollary}

We remark that the first part of this corollary is proven
in~\cite{Shi:01} for the special case where $\CycleSet(Y)$ is a
circle. The second part is well-known (see, \eg~\cite{Humphreys:92}). 

Let $P = (v_1,v_2,\ldots,v_k)$ be a (possibly closed) simple path in
$Y$. The map
\begin{equation}
\label{eq:nu}
\nu_P \colon \Acyc(Y) \longrightarrow \Z
\end{equation}
evaluated at $O_Y$ is the number of edges of the form
$\{v_i,v_{i+1}\}$ in $Y$ oriented as $(v_i, v_{i+1})$ in $O_Y$
(positive edges) minus the number of edges oriented as $(v_{i+1},v_i)$
in $O_Y$ (negative edges).
\begin{lemma}
\label{lem:nu}
Let $P$ be a simple closed path in the undirected graph $Y$. The map
$\nu_P$ extends to a map $\nu^*_P \colon \Acyc(Y)/\simkappa
\longrightarrow \Z$. 
\end{lemma}
\begin{proof}
Let $c(O_Y) = O'_Y$ where $c = c_v$ is a click of a single vertex $v$. If
$v$ is not an element of $P$ then clearly $\nu_P(O_Y) =
\nu_P(O'_Y)$. On the other hand, if $v$ is contained in $P$ then $c$
maps one positive edge into a negative edge and vice versa. The
general case follows by induction on the length of the click sequence.
\qed
\end{proof}
Lemma~\ref{lem:nu} will be used extensively in the proof of the main
result in the next section. 


\section{Proof of the Main Theorem}

From Proposition~\ref{prop:non-cycle-edge} it is clear that for the
computation of $\kappa(Y)$ all bridges may be omitted. We now turn our
attention to the role played by cycle-edges in determining
$\kappa(Y)$ and to the proof of the recursion relation
\begin{equation}
  \kappa(Y) = \kappa(Y_e')+\kappa(Y_e'')
\end{equation}
of Theorem~\ref{thm:collapse} valid for any cycle-edge $e$ of $Y$.  We
set $e = \{v,w\}$ in the following.

First, define $\iota_1 \colon \Acyc(Y_e'') \longrightarrow \Acyc(Y)$
as the map that sends $O_{Y''}\in\Acyc(Y_e'')$ to $O_Y\in\Acyc(Y)$ for
which $O_Y(e) = (v,w)$ and for which all other edge orientations are
inherited. The map $\iota_2 \colon \Acyc(Y_e'') \longrightarrow
\Acyc(Y)$ is defined analogously, but orients $e$ as $(w,v)$.
\begin{proposition}
The maps $\iota_{1,2} \colon \Acyc(Y_e'') \longrightarrow \Acyc(Y)$
extend to well-defined maps
\begin{equation}
\label{eq:iota-extend}
\iota^*_{1,2} \colon \Acyc(Y_e'')/\simkappa 
   \longrightarrow \Acyc(Y)/\simkappa \;.
\end{equation}
\end{proposition}
\begin{proof}
For $\phi \in \{\iota_1, \iota_2\}$ and for any click-sequence $c$
of $O_{Y''} \in \Acyc(Y'')$ we have the commutative diagram
\begin{equation}
\label{eq:comm}
\xymatrix{
O_{Y''} \ar[rr]^{c} \ar[d]_\phi && O_{Y''}'\ar[d]^{\phi}  \cr
O_{Y} \ar[rr]^{c'} && O_{Y}'  
}
\end{equation}
where the click-sequence $c'$ over $Y$ is constructed from the
click-sequence $c$ over $Y''$ by insertion of $w$ after
(resp. before) every occurrence of $v$ in $c$ for $\iota_1$
(resp. $\iota_2$). 
\qed
\end{proof}

\begin{proposition}
\label{prop:iota-diff}
Let $e$ be a cycle-edge.  For any $[O_{Y''}] \in \Acyc(Y'')$ we have
$\iota^*_1([O_{Y''}]) \ne \iota^*_2([O_{Y''}])$. 
\end{proposition}
\begin{proof}
Let $P$ be any simple closed path containing $e$ and oriented so as to
include $(v,w)$. From the definition of $\iota_1$ and $\iota_2$ we
conclude that $\nu_P\bigl(\iota_1(O_{Y''})\bigr) =
\nu_P(\iota_2\bigl(O_{Y''})\bigr) + 2$, and the proposition follows by
Lemma~\ref{lem:nu}.
\end{proof}

\begin{proposition}
The maps $\iota^*_{1,2}$ are injections. 
\end{proposition}
\begin{proof}
We prove the statement for $\iota^*_1$. The proof for $\iota^*_2$ is
analogous.  Assume $[O_{Y''}] \not\sim_\kappa [O'_{Y''}]$
both map to $[O_Y]$ under $\iota^*_1$. By construction, any elements
$O_{Y''}$ and $O_{Y''}'$ of the respective $\kappa$-classes have
$\iota_1$-images with $e$ oriented as $(v,w)$. Moreover, for any image
point of $\iota_1$ there is no directed path from $v$ to $w$ of length
$\ge 2$, and there is no directed path from $w$ to $v$. We may also
assume that $v$ is a source in both $O_{Y''}$ and $O'_{Y''}$\;. From
this it is clear that $v$ and $w$ belong to successive layers in the
acyclic orientations $O_Y$ and $O_Y'$.

Let $c$ be a click-sequence taking $\iota_1(O_{Y''})$ to
$\iota_1(O_{Y''}')$. Again by construction, we may assume that any
occurrence of $v$ in $c$ is immediately followed by $w$. This follows
since $v$ and $w$ have to occur equally many times in $c$, and from
the fact that $v$ and $w$ belong to successive layers in $O_Y$ and
$O_Y'$. If $v$ and $w$ were not consecutive in $c$ it could only be
because $c$ is of the form $c = \ldots v\, v_1  \ldots v_k\, w_1 \ldots
w_r \, w \ldots$ where the $v_i$'s belong to the same layer as $v$ and the
$w_i$'s belong to the same layer as $w$. A layer is in particular an
independent set, and it is clear that the sequence $c' = \ldots v_1
\ldots v_k\, v\, w\, w_1 \ldots w_r \ldots$ obtained from $c$ by changing
the order of $v$ and $w$ also maps $O_Y$ to $O_Y'$.

The click-sequence $c''$ obtained from $c'$ by deleting every
occurrence of $w$ is a click-sequence mapping $O_{Y''}$ to $O_{Y''}'$,
which contradicts the assumption that $[O_{Y''}] \not\sim_\kappa
[O_{Y''}]'$.
\qed
\end{proof}

Consequently, any $\kappa$-class $[O_Y]$ contains at most one set of
the form $\iota_1([O_{Y''}])$, and at most one set of the form
$\iota_2([O_{Y''}])$. 

\begin{proposition}
\label{prop:c-simple}
Let $e$ be a cycle-edge of the undirected graph $Y$.  For each pair of
distinct $\kappa$-classes $[O_Y]$ and $[O_Y]'$ there is at most one
$\kappa$-class $[O_{Y''}]$ such that $\{\iota^*_1([O_{Y''}]),
\iota^*_2([O_{Y''}]) \} = \{[O_Y], [O_Y]'\} $.
\end{proposition}
\begin{proof}
Assume this is not not the case, and that there in fact is another
class $[O_{Y''}]'$ with the same property. Since both maps
$\iota^*_{1,2}$ are injective it then follows (up to relabeling) that
$\iota^*_1([O_{Y''}]) = \iota^*_2([O_{Y''}]') = [O_Y]$ and 
$\iota^*_1([O_{Y''}]') = \iota^*_2([O_{Y''}]) = [O_Y]'$. 
By the same argument as in the proof of
Proposition~\ref{prop:iota-diff} it follows using $[O_{Y''}]$ that
$\nu^*_P([O_Y]) = \nu^*_P([O_Y]') + 2$. On the other hand, by using
$[O_{Y''}]'$ if follows that that $\nu^*_P([O_Y]') = \nu^*_P([O_Y]) +
2$, which is impossible.
\qed
\end{proof}

\begin{definition}
\label{def:collapse}
Let $e$ be a cycle-edge of the undirected graph $Y$. The collapse
graph $\mathfrak{C}_e(Y)$ of $Y$ and $e$ is the graph with vertex set
$\Acyc(Y)/\simkappa$ and edge set $\bigl\{\{\iota^*_1([O_{Y''}]),
\iota^*_2([O_{Y''}]) \} \,\bigl|\, [O_{Y''}] \in \Acyc(Y'')/\simkappa
\bigr\}$.
\end{definition}

Note that by Proposition~\ref{prop:c-simple}, the graph
$\mathfrak{C}_e(Y)$ is simple, and by
Proposition~\ref{prop:iota-diff}, it has no singleton edges (\ie
self-loops).

A \emph{line graph} on $k$ vertices has vertex set $\{1,2,\ldots,k\}$
and edges $\{i,i+1\}$ for $1\le i \le k-1$.

\begin{proposition}
Let $e$ be a cycle-edge of the undirected graph $Y$. The collapse
graph $\mathfrak{C}_e(Y)$ is isomorphic to a disjoint collection of line
graphs.
\end{proposition}
\begin{proof}
By Definition~\ref{def:collapse} and the remark following it, each
vertex in the collapse graph has degree $\le 2$. Thus it is sufficient
to show that $\mathfrak{C}_e(Y)$ contains no cycles.
By the now familiar argument using $\nu_P$ for some path containing
$e$, two adjacent $\kappa$-classes in $\mathfrak{C}_e(Y)$ differ in
their $\nu^*$-values by precisely $2$. Assume $\mathfrak{C}_e(Y)$
contains the subgraph (up to relabeling)
\begin{equation}
\label{eq:line}
\xymatrix@R5pt{
 A'' \ar@{-}[rr] && A \ar@{-}[rr] && A'  \cr
\iota^*_1(A_1) && \iota^*_2(A_1) = \iota^*_1(A_2) && \iota^*_2(A_2) 
} \;,
\end{equation}
for unique $\kappa$-classes $A_{1,2}\in\Acyc(Y'')/\simkappa$. 
Clearly $\nu^*_P(A'') = \nu^*_P(A) + 2$ and $\nu^*_P(A') = \nu^*_P(A)
- 2$.  We now have the following situation: $(i)$ the $\nu^*$-values of
adjacent vertices in $\mathfrak{C}_e(Y)$ differ by precisely $2$, and
$(ii)$ the value of $\nu*$ increases by $2$ across each edge in the
$A''$-direction relative to $A$ and decreases by $2$ across each edge
in the $A'$-direction relative to $A$. If the subgraph
in~\eqref{eq:line} was a part of cycle of length $\ge 3$ in
$\mathfrak{C}_e(Y)$, then by $(ii)$ there must be some pair of
adjacent vertices for which $\nu^*$ differs by at least $4$. But this is
impossible by $(i)$, hence $\mathfrak{C}_e(Y)$ contains no cycles and
the proof is complete.
\qed
\end{proof}

\begin{proposition}
\label{prop:compstruct}
Let $e$ be a cycle-edge of the undirected graph $Y$. Then
$\kappa$-classes on the same connected component in
$\mathfrak{C}_e(Y)$ are contained in the same $\kappa$-class in
$\Acyc(Y')/\simkappa$.
\end{proposition}

\begin{proof}
It is sufficient to show this for adjacent vertices in
$\mathfrak{C}_e(Y)$ -- the general result follows by induction.
Clearly, $O_Y \simkappa O_Y'$ implies $O_{Y'} \simkappa
O_{Y'}'$. Adjacent vertices in $\mathfrak{C}_e(Y)$ contain elements
that only differ in their orientations of $e$ and hence become
$\kappa$-equivalent upon deletion of $e$. The proof follows. 
\qed
\end{proof}

\begin{proposition}
\label{prop:yprime}
There is a bijection between the connected components in
$\mathfrak{C}_e(Y)$ and $\Acyc(Y_e')/\simkappa$.
\end{proposition}

\begin{proof}
  Let $n_c$ denote the number of connected components of
  $\mathfrak{C}_e(Y)$. By Proposition~\ref{prop:compstruct} if
  $[O_Y]$ and $[O_Y']$ are connected in $\mathfrak{C}_e(Y)$ then both
  classes are contained in the same $\kappa$-class over $Y'$, and thus
  $n_c \le \kappa(Y')$.

  It is clear that a $\kappa$-class contains all acyclic orientations
  for which there are representative permutations that are related by
  a sequence of adjacent transpositions of non-connected vertices in
  $Y$ and cyclic shifts.  Upon deletion of the cycle-edge $e$ the
  adjacent transposition of the endpoints of $e$ becomes permissible,
  and thus two distinct $\kappa$-classes in $Y$ containing acyclic
  orientations that only differ on $e$ are contained within the same
  $\kappa$-class over $Y'$.  By reference to the underlying
  permutations, it follows that two $\kappa$-classes in $Y$ are
  contained within the same $\kappa$-class in $Y'$ if and only if
  there is a sequence of $\kappa$-classes in $Y$ where consecutive
  elements in the sequence contain acyclic orientations that differ
  precisely on $e$.  By the definition of $\mathfrak{C}_e(Y)$ it
  follows that all $\kappa$-classes over $Y$ that merge to be
  contained within one $\kappa$-class in $Y'$ upon deletion of $e$ are
  contained within the same connected component of
  $\mathfrak{C}_e(Y)$, and thus $n_c \ge \kappa(Y')$.
  \qed
\end{proof}

\begin{proof}[Theorem~\ref{thm:collapse}]
  Upon deletion of $e$ in $Y$ two or more $\kappa$-classes over $Y$
  may merge to be contained with the same $\kappa$-class over $Y'$. By
  Proposition~\ref{prop:yprime} the number of $\kappa$-classes over
  $Y'$ equals the number of connected components in
  $\mathfrak{C}_e(Y)$. If a connected component in $\mathfrak{C}_e(Y)$
  contains $m$ distinct $\kappa$-classes of $Y$ then by
  Proposition~\ref{prop:c-simple} there are $m-1$ unique corresponding
  $\kappa$-classes over $Y''$. Thus for each component of
  $\mathfrak{C}_e(Y)$ we have a relation precisely of the
  form~\eqref{eq:recursion} for the $\kappa$-classes involved. The
  theorem now follows.
\qed
\end{proof}


\section{Related Enumeration Problems}
\label{sec:tutte}

In this section we relate the problem of computing $\kappa(Y) =
\card{\Acyc(Y)/\simkappa}$ to two other enumeration problems where the
same recurrence holds. We will show how these problems are equivalent,
and additionally, how they all can be computed through an evaluation
of the Tutte polynomial. As a corollary we obtain a transversal of
$\Acyc(Y)/\simkappa$.

In~\cite{BChen:07} the notion of \emph{cut-equivalence} of acyclic
orientations is studied. Recall that a \emph{cut} of a graph $Y$ is a
partition of the vertex set into two classes $\vset[Y]=V_1\sqcup V_2$,
and where $[V_1,V_2]$ is the set of (cut-)edges between $V_1$ and
$V_2$.  A cut of a directed graph of the form $G(O_Y)$, which we
simply refer to as a cut of $O_Y$, is \emph{oriented} with respect to
$O_Y$ if the edges of $[V_1, V_2]$ are all directed from $V_1$ to
$V_2$, or are all directed from $V_2$ to $V_1$.
\begin{definition}
  Two acyclic orientations $O_Y$ and $O'_Y$ are \emph{cut-equivalent}
  if the set $\{e\in\vset[Y] \mid O_Y(e) \ne O_Y'(e) \}$ is $(i)$ empty
  or is $(ii)$ an oriented cut with respect to either $O_Y$ or $O_Y'$.
\end{definition}

The study of cut-equivalence in~\cite{BChen:07} was done outside the
setting of Coxeter theory, and here we provide the connection.

\begin{proposition}
  Two acyclic orientations of $Y$ are $\kappa$-equivalent if and only
  if they are cut-equivalent. 
\end{proposition}

\begin{proof}
  Suppose distinct elements $O_Y$ and $O'_Y$ in $\Acyc(Y)$ are
  cut-equivalent, and without loss of generality, that all edges of
  $[V_1,V_2]$ are oriented from $V_1$ to $V_2$ in $O_Y$.  A
  click-sequence containing each vertex of $V_1$ precisely once maps
  $O_Y$ to $O_Y'$, thus $O_Y\simkappa O'_Y$.

  Conversely, suppose that $O_Y\simkappa O'_Y$, where $O'_Y$ is
  obtained from $O_Y$ by a click-sequence containing a single vertex
  $v$. Then $O_Y$ and $O'_Y$ are cut-equivalent, with the partition
  being $\{v\}\sqcup\vset[Y]\setminus\{v\}$. 
  \qed
\end{proof}

Obviously, the recurrence relation in~\eqref{eq:recursion} holds for
the enumeration of both cut-equivalence and $\kappa$-equivalence
classes.

\begin{definition}
  The \emph{Tutte polynomial} of an undirected graph $Y$ is defined as
  follows. If $Y$ has $b$ bridges, $\ell$ loops, and no cycle-edges,
  then $T(Y,x,y)=x^by^\ell$. If $e$ is a cycle-edge of $Y$, then
  \begin{equation*}
    T(Y,x,y)=T(Y'_e,x,y)+T(Y''_e,x,y)\;.
  \end{equation*}
\end{definition}

We remark that it is well-known that the number of acyclic
orientations of a graph $Y$ can be evaluated as $\alpha(Y) =
T(Y,2,0)$. It was shown in~\cite{BChen:07} that the number of
cut-equivalence classes can be computed through an evaluation of the
Tutte polynomial as $T(Y,1,0)$, and thus $\kappa(Y) = T(Y,1,0)$. 

It is known that $T(Y,1,0)$ counts several other quantities, some of
which can be found in~\cite{Gioan:07}. One of these is
$\card{\Acyc_v(Y)}$, the number of acyclic orientations of $Y$ where a
fixed vertex $v$ is the unique source. As the next result shows, there
is a bijection between $\Acyc(Y)/\simkappa$ and $\Acyc_v(Y)$.

\begin{proposition}
  \label{prop:uniquesource}
  Let $Y$ be a connected graph $Y$. For any fixed $v\in\vset[Y]$ there
  is a bijection
\begin{equation*}
  \phi_v \colon \Acyc(Y)/\simkappa \longrightarrow \Acyc_v(Y)\;.
\end{equation*}
\end{proposition}

\begin{proof}
  Since $\kappa(Y) = \card{\Acyc(Y)/\simkappa} = T(Y,1,0) =
  \card{\Acyc_v(Y)}$ it is sufficient to show that that there is a
  surjection $\beta_v \colon \Acyc(Y)/\simkappa \longrightarrow
  \Acyc_v(Y)$. Let $A \in \Acyc(Y)/\simkappa$, let $O_Y \in A$, and
  let $c$ be a maximal length click-sequence not containing the vertex
  $v$. This click-sequence is finite since $Y$ is connected. The
  orientation $c(O_Y)$ has $v$ as a unique source, since otherwise $c$
  would not be maximal, and the proof is complete. \qed
\end{proof}

From Proposition~\ref{prop:uniquesource} we immediately obtain:

\begin{corollary}
  For any vertex $v$ of $Y$ the set $\Acyc_v(Y)$ is a transversal of
  $\Acyc(Y)/\simkappa$.
\end{corollary}

In light of the results in this section, the recurrence
in~\eqref{eq:recursion} may also be proven by showing that the map
$\beta_v$ is injective. However, our proof offers insight into the
structure of the $\kappa$-classes, and it is our hope that this may
lead to new techniques for studying conjugacy classes of Coxeter
groups.


\begin{acknowledgement}
The first author would like to thank Jon McCammond for many helpful
discussions. Both authors are grateful to the NDSSL group at Virginia
Tech for the support of this research. Special thanks to Vic Reiner
and William~Y.~C. Chen for valuable advise regarding the content,
preparation, and structure of this paper. The work was partially
supported by Fields Institute in Toronto, Canada.
\end{acknowledgement}


\end{document}